\input amstex
\loadbold
\documentstyle{amsppt}

\input xy  %try \input xypic if that doesn't work
\xyoption{all}

\define\mathbb{\Bbb}
\define\mathscr{\Cal}
\define\mathcal{\Cal}
\define\mathfrak{\frak}
\define\em{\it}

\NoBlackBoxes

\topmatter
\title 
Eigenvalues of majorized Hermitian matrices and Littlewood-Richardson 
coefficients.
\endtitle
\rightheadtext{Eigenvalues of majorized Hermitian matrices}
\author William Fulton \endauthor
\address  University of Michigan,  
	Ann Arbor, MI 48109-1109  \endaddress
\email wfulton\@math.lsa.umich.edu \endemail
\date September 21, 2000 \enddate
\thanks The author was partly supported by NSF Grant \#DMS9970435 
\endthanks

\abstract
Answering a question raised by S\. Friedland, we show that the possible 
eigenvalues of Hermitian matrices (or compact operators) $A$, 
$B$, and $C$ with $C \leq A + B$ are given by the same inequalities as in 
Klyachko's theorem for the case where $C = A + B$, except that the equality 
corresponding to $\operatorname{tr}(C) = \operatorname{tr}(A) + 
\operatorname{tr}(B)$ is replaced by the inequality corresponding to 
$\operatorname{tr}(C) \leq \operatorname{tr}(A) + \operatorname{tr}(B)$.  
The possible types of finitely generated torsion modules $A$, $B$, and $C$ 
over a discrete valuation ring such that there is an exact sequence $B \to C 
\to A$ are characterized by the same inequalities. 
\endabstract

%\toc
%\widestnumber\head{2}	
%\head   1.  Prelimaries from Schubert calculus  \endhead
%\head   2.  Theorems and their proofs \endhead
%\endtoc

\endtopmatter
\document

\head Introduction  \endhead
A\. Klyachko [10], cf. [13], [6], [7], has shown that the possible 
eigenvalues $\alpha$, $\beta$, $\gamma$ of Hermitian $n$ by $n$ 
matrices $A$, $B$, $C$ with $C = A + B$ are characterized by a certain 
list of inequalities, together with the trace equality $\sum \gamma_i = 
\sum \alpha_i + \sum \beta_i$.
Combined with the solution of the saturation conjecture by A\. Knutson and T\. 
Tao [11], cf. [2], [3], one can show that this list of inequalities is 
exactly that conjectured by A\. Horn [8].  P\. Belkale [1]  
showed that that a smaller list of inequalities is sufficient, and 
Knutson, Tao, and C\. Woodward [12] have announced that this smaller list 
is actually minimal.  This story, with more of 
the history and references, can be found in [7].

Eigenvalues (and all n-tuples of real numbers) are always written in 
descending order, so we start with the $3n-3$ inequalities
$$
\alpha_1 \geq \alpha_2 \geq \ldots \geq \alpha_n, \quad
\beta_1 \geq \beta_2 \geq \ldots \geq \beta_n, \quad
\gamma_1 \geq \gamma_2 \geq \ldots \geq \gamma_n.   \tag{$\dagger$}
$$
S\. Friedland [4] subsequently gave a characterization of the inequalities 
satisfied by the eigenvalues $\alpha$, $\beta$, $\gamma$ of Hermitian $n$ 
by $n$ matrices $A$, $B$, and $C$ with $C \leq A + B$, i.e., such that 
$A+B-C$ is positive semidefinite.  His inequalities included those of 
Klyachko, together with the trace inequality 
 $$
\sum_{i=1}^n \gamma_i \, \leq \,\sum_{i=1}^n \alpha_i + \sum_{i=1}^n 
\beta_i,    \tag{$\dagger \dagger$}
$$
but they also included some other inequalities that are less easy to describe.  
For $n=2$ and $n=3$, however, he showed that these extra inequalities are 
superfluous.  

In this paper we show that these extra inequalities are always superfluous, so 
that the natural generalization of Klyachko et al\. 
to the majorization situation 
is valid.

The inequalities of Horn or Klyachko have the form 
$$
\sum_{k \in K} \gamma_k \, \leq \, \sum_{i \in I} \alpha_i + \sum_{i \in J}  
\beta_j \tag{$\dagger_{IJK}$}
$$
where $(I,J,K)$ vary over certain subsets of the same cardinality $r$ of the 
index set $[n] = \{1, \ldots ,n\}$, $1 \leq r \leq n-1$.  If, for $I = \{i_1 < 
\ldots < i_r\}$ we set 
$$
\lambda(I) = (i_r - r, \ldots, i_2 - 2, i_1 - 1),
$$
then $(I,J,K)$ is on the Horn/Klyachko list exactly when the 
Littlewood-Richardson
coefficient $c_{\lambda(I)\,\lambda(J)}^{\,\,\,\,\lambda(K)}$ of the 
three corresponding partitions is positive.  The minimal list of inequalities 
of Belkale-Knutson-Tao-Woodward consists of those ($\dagger_{IJK}$) for 
which $c_{\lambda(I)\,\lambda(J)}^{\,\,\,\,\lambda(K)} = 1$, together with the 
trace equality and the inequalities ($\dagger$), for $n \geq 3$.  (For $n=2$ 
the inequalites ($\dagger$) follow from the others.)

\proclaim{Theorem 1}  A triple $\alpha$, $\beta$, $\gamma$ satisfying 
$(\dagger)$ occurs as the eigenvalues of $n$ by $n$ Hermitian matrices 
$A$, $B$, $C$ with $C \leq A + B$ if and only if they satisfy $(\dagger 
\dagger)$ and $(\dagger_{IJK})$ for all $(I,J,K)$ of cardinality $r < n$ 
such that $c_{\lambda(I)\,\lambda(J)}^{\,\,\,\,\lambda(K)} = 1$. \endproclaim

In fact, granting the minimality assertion of [12], we show that these 
inequalities ($\dagger$), ($\dagger \dagger$), and ($\dagger_{IJK}$) are a 
minimal list of inequalities for this majorization problem, for all $n \geq 1$.  

As in [4], this generalizes to give an infinite list of inequalities for 
eigenvalues of compact self-adjoint nonnegative operators $A$, $B$, $C$ 
with $C \leq A + B$.  

As in [6] and [7], the same assertions are valid for real symmetric  
matrices (or operators) or for quaternionic Hermitian matrices.  All the 
results extend to sums of $m$ rather than two factors; these generalizations 
are given in Section 2.

As explained in [7], Klyachko's theorem has analogues in several other 
areas of mathematics, such as representation theory, algebra, and 
combinatorics.  The same is true of the stronger theorems presented here.  
In particular, the inequalities in Theorem 1, when restricted to the case 
where $\alpha$, $\beta$, and $\gamma$ are partitions, are equivalent 
to the conditions that $\gamma$ is contained in 
\footnote{A partition $\widetilde{\gamma}$ {\em contains} a partition $\gamma$ 
if $\widetilde{\gamma}_i \geq \gamma_i$ for all $i$, i.e., the Young 
diagram of $\widetilde{\gamma}$ contains that of $\gamma$.  In this case 
we write $\widetilde{\gamma} \supset \gamma$.} 
some partition $\widetilde{\gamma}$ such that the Littlewood-
Richardson coefficient $c_{\alpha\,\beta}^{\,\,\widetilde{\gamma}}$ is 
positive, and also to the condition that $\alpha$ and $\beta$ contain 
partitions $\widetilde{\alpha}$ and $\widetilde{\beta}$  such that 
$c_{\widetilde{\alpha}\,\widetilde{\beta}}^{\,\,\gamma}$ is positive.  A 
consequence is a simple criterion for there to be 
homomorphisms putting three given finite abelian $p$-groups in an 
exact sequence.  More generally, recall that the {\em type} of a finitely 
generated torsion module $A$ over any discrete valuation ring $R$ is the 
partition $\alpha \, : \, \alpha_1 \geq \ldots \geq \alpha_n$ 
such that $A$ is 
isomorphic to $R/\mathfrak{p}^{\alpha_1} \oplus \ldots \oplus
R/\mathfrak{p}^{\alpha_n}$, where $\mathfrak{p}$ is the maximal ideal of 
$R$.  Then we have:

\proclaim{Theorem 2}  For any discrete valuation ring $R$, there is an 
exact sequence  
$$
B \to C \to A 
$$
of $R$-modules of types $\beta$, $\gamma$, and $\alpha$ (all partitions of 
lengths at most $n$) if and only if they satisfy the inequalities $(\dagger)$, 
$(\dagger \dagger)$, and $(\dagger_{IJK})$ for all $(I,J,K)$ of cardinality $r 
< n$ such that $c_{\lambda(I)\,\lambda(J)}^{\,\,\,\,\lambda(K)} = 1$. 
\endproclaim

This is a minimal list of inequalities.  More general results appear in 
Section 2.

One naturally expects Theorem 1 
to be a consequence of Klyachko's theorem, since 
an inequality $C \leq A(1) + \ldots + A(m)$ is equivalent to the existence of an 
$A(m+1)$ with nonpositive eigenvalues such that $C = A(1) + \ldots + 
A(m+1)$.  This means that the polyhedral cone describing the majorization 
problem is a projection of a polyhedral cone describing the equality 
problem with one more factor.  Computing explicit inequalities to describe 
the projection of a polyhedral cone, however, is seldom easy, and it is 
exactly this that leads to the extraneous and inexplicit inequalities of [4].  

Our proof does depend on the Klyachko case of equality, but not by a 
projection.  The essential point is to show that if some inequality 
($\dagger_{IJK}$) is an equality, then the triple $(\alpha,\beta,\gamma)$ 
splits into two triples, one $(\alpha',\beta',\gamma')$ of length $r$ 
(consisting of those $\alpha_i$ for $i \in I$, $\beta_j$ for $j \in J$ and 
$\gamma_k$ for $k \in K$), and one $(\alpha'', \beta'', \gamma'')$ of 
length $n-r$ (consisting of the others).  We show that 
$(\alpha',\beta',\gamma')$ satisfies the conditions to be the eigenvalues of 
$r$ by $r$ Hermitian matrices $A'$, $B'$, $C'$ with $C' = A' + B'$, and, by 
induction, that $(\alpha'',\beta'',\gamma'')$ satisfies the conditions to be 
the eigenvalues of $n-r$ by $n-r$ Hermitian matrices $A''$, $B''$, $C''$ with 
$C'' \leq A'' + B''$.  The direct sums $A = A' \oplus A''$, $B = B' \oplus B''$, 
and $C = C' \oplus C''$ then do the trick.  

To prove that these inductive conditions are satisfied involves a little 
Schubert calculus, which is carried out in the first section.  This deduces 
some nonzero intersections in a Grassmann variety 
$\operatorname{Gr}(r,n)$ of $r$-planes in $\mathbb{C}^n$ from nonzero 
intersections in $\operatorname{Gr}(p,r)$ or $\operatorname{Gr}(p,n-r)$ 
for smaller $p$.  This can be regarded as a contribution from the Schubert 
calculus side of the general problem, emphasized in [7], of 
understanding why all the conditions in the variations of Klyachko's theorem 
are inductive, i.e., their answers for given $n$ are determined by knowing 
the answers to the same questions for smaller $n$.

We thank H. Derksen, J. Weyman, J. Harris, and A. Buch for stimulating 
discussions that led to the proof given here, and S. Friedland for raising the 
question and suggesting improvements.

\head Section 1.  Preliminaries from Schubert calculus \endhead

We start by fixing some notation.

A subset $I$ of $[n] = \{1, \ldots ,n\}$ is always written in increasing order, 
so that  
$I = \{i_1 < i_2 < \ldots < i_r\}$.  For a subset $P$ of $[r]$ of cardinality 
$x$, set
$$
I_P = \{ i_p \mid p \in P \} = \{ i_{p_1} < \ldots < i_{p_x} \},
$$
a subset of $[n]$ of cardinality $x$.  For $P$ a subset of $[n-r]$ of cardinality 
$y$, set 
$$
I_P^{+} = I \cup (I^c)_P,
$$
where $I^c$ denotes the complement of $I$ in $[n]$; this is a subset of $[n]$ 
of cardinality $r+y$.  For $m$-tuples $\mathscr{I} = (I(1), \ldots , I(m))$ 
and $\mathscr{P} = (P(1), \ldots , P(m))$ of such subsets, we write 
$\mathscr{I}_{\mathscr{P}}$ for $(I(1)_{P(1)}, \ldots , I(m)_{P(m)})$, and  
$\mathscr{I}_{\mathscr{P}}^{+}$ for $(I(1)_{P(1)}^{+}, \ldots , 
I(m)_{P(m)}^{+})$.

For subsets $H$ and $I$ of $[n]$ of the same cardinality $r$, we write $H 
\leq I$ if $H = \{ h_1 < \ldots  < h_r \}$ 
and $I = \{ i_1 < \ldots < i_r \}$ with 
$h_a \leq i_a$ for $1 \leq a \leq r$.  Equivalently, 
$
|H \cap [k]| \geq |I \cap [k]|$ for $1 \leq k \leq n$.

Let $V$ be an $n$-dimensional complex vector space, and let
$$
F_{\sssize{\bullet}} \,\, : \,\, 0 = F_0 \subset F_1 \subset \ldots \subset 
F_n = V 
$$
be a complete flag of subspaces of $V$.  For a subset $I$ of $[n]$ of 
cardinality $r$, the {\em Schubert variety}  
$\Omega_I(F_{\sssize{\bullet}})$ is the subvariety of the Grassmannian 
$\operatorname{Gr}(r,V)$ of $r$-planes in $V$ defined by
$$
\Omega_I(F_{\sssize{\bullet}}) = \{ L \in \operatorname{Gr}(r,V) \mid 
\operatorname{dim}(L \cap F_{i_k}) \geq k \,\,\, {\text {for }} 1 \leq k \leq 
r \}.
$$
The class of $\Omega_I(F_{\sssize{\bullet}})$ in the cohomology ring 
$H^*(\operatorname{Gr}(r,n)) = H^*(\operatorname{Gr}(r,V))$, which is 
independent of choice of the flag $F_{\sssize{\bullet}}$, is denoted 
$\omega_I$.  If $\lambda = (\lambda_1, \ldots, \lambda_r)$, 
with $\lambda_k = n-r+k-i_k$, this class 
is also denoted $\sigma_\lambda$; the complex codimension of 
$\Omega_I(F_{\sssize{\bullet}})$ is $|\lambda| = \sum_{i=1}^r 
\lambda_i$, so we have 
$$
\sigma_\lambda =  \omega_I = 
[\Omega_I(F_{\sssize{\bullet}})] \in 
H^{2|\lambda|}(\operatorname{Gr}(r,n)).
$$
The Schubert variety $\Omega_I(F_{\sssize{\bullet}})$ is the closure of the 
{\em Schubert cell}  
$\Omega_I^{\circ}(F_{\sssize{\bullet}})$, which is the set of subspaces $L$ 
such that the jumps in the sequence
$$
0 = \operatorname{dim}(L \cap F_0) \leq \operatorname{dim}(L \cap F_1) 
\leq \ldots \leq \operatorname{dim}(L \cap F_n) = r
$$ 
occur exactly at the integers in $I$; that is, $I = \{ i \in [n] \mid L\cap F_i 
\neq L \cap F_{i-1} \}$.  For a fixed flag, every $r$-dimensional subspace 
$L$ of $V$ belongs to a unique Schubert cell.  We use the well-known fact 
(cf. [5], \S 10) that 
$$
\Omega_I(F_{\sssize{\bullet}}) = \bigcup_{J \leq I} 
\Omega_J^{\circ}(F_{\sssize{\bullet}}) \, .
$$

\proclaim{Lemma 1}  Let $H$ and $I$ be subsets of $[n]$ of cardinality $r$ 
with $H \leq I$.  
\roster \widestnumber\item{(ii)}
\item"{(i)}" If $P$ and $Q$ are subsets of $[r]$ of cardinality $x$ with $P 
\leq Q$, then $H_P \leq I_Q$ (as subsets of $[n]$ of cardinality $x$).
\item"{(ii)}" If $P$ and $Q$ are subsets of $[n-r]$ of cardinality $y$ with $P 
\leq Q$, then $H_P^{+} \leq I_Q^{+}$ (as subsets of $[n]$ of cardinality 
$r+y$). 
\endroster
\endproclaim
\demo{Proof}  The proof of (i) is trivial: $h_{p_k} \leq h_{q_k} \leq 
i_{q_k}$.  With the assumptions of (ii), it follows from (i) that $I_P^{+} \leq 
I_Q^{+}$, so it suffices to show that $H_P^{+} \leq I_P^{+}$.  Any $H \leq 
I$ can be obtained from $I$ by a succession of moves, each of which 
decreases one integer, leaving the others alone.  We may therefore 
assume that, for some $k \in [r]$ and some $m \in [n]$, $h_a = i_a$ for $a 
\neq k$, and $h_k = m-1$, $i_k = m$.  The complementary sequences 
satisfy $h_a^c = i_a^c$ for $a \neq m-k$, while $i_{m-k}^c = m-1$ and 
$h_{m-k}^c = m$.  Therefore $H_P^{+} = I_P^{+}$ if $P$ contains $m-k$, 
and $H_P^{+} < I_P^{+}$ otherwise.  \enddemo

Let $U$ be a subspace of $V$ of dimension $r$.  A complete flag 
$F_{\sssize{\bullet}}$ on $V$ determines a complete flag 
$F_{\sssize{\bullet}} U$ on $U$ and a complete flag $F_{\sssize{\bullet}} 
W$ on $W = V/U$, which are defined as follows.  Let $I = \{ i_1 < \ldots 
< i_r \}$ 
be the set such that 
 $U$ is in $\Omega_I^\circ(F_{\sssize{\bullet}})$.  Set 
$$
F_kU = F_{i_k} \cap U, \quad 1 \leq k \leq r.
$$
With $I^c = \{ i_1^c < \ldots < i_{n-r}^c \}$, set
 $$
F_kW = (U + F_{i_k^c}) / U, \quad 1 \leq k \leq n-r.
$$
From the isomorphism $(U + F_i)/U \cong F_i/F_i\cap U$ it follows that 
$F_i \cap U = F_{i-1} \cap U$ if and only if $U + F_i \neq U + F_{i-1}$.

\proclaim{Lemma 2}  Let $U$ be in a Schubert variety 
$\Omega_I(F_{\sssize{\bullet}})$ in $\operatorname{Gr}(r,V)$. 
\roster\widestnumber\item{(ii)}
\item"{(i)}"  If $X$ is a subspace of $U$ of dimension $x$, with $X$ in 
$\Omega_P(F_{\sssize{\bullet}} U)$ in $\operatorname{Gr}(x,U)$, then 
$X$ is in $\Omega_{I_P}(F_{\sssize{\bullet}})$ in $\operatorname{Gr}(x,V)$.
\item"{(ii)}" If $Y$ is a subspace of $W$ of dimension $y$, with $Y$ in 
$\Omega_P(F_{\sssize{\bullet}} W)$ in $\operatorname{Gr}(y,W)$, and $Y = 
Z/U$, then $Z$ is in $\Omega_{I_P^{+}}(F_{\sssize{\bullet}})$ in 
$\operatorname{Gr}(r+y,V)$.
\endroster
\endproclaim
\demo{Proof}  For (i), if $U \in \Omega_H^\circ(F_{\sssize{\bullet}})$ and 
$X \in \Omega_Q^\circ(F_{\sssize{\bullet}} U)$, then by Lemma 1(i), $H_Q 
\leq I_P$, so $\Omega_{H_Q}(F_{\sssize{\bullet}}) \subset 
\Omega_{I_P}(F_{\sssize{\bullet}})$.  So we may assume 
$U \in \Omega_I^\circ(F_{\sssize{\bullet}})$ and $X \in 
\Omega_P^\circ(F_{\sssize{\bullet}} U)$.  Then the jumps in the sequence 
$\operatorname{dim}(U \cap F_i)$, $1 \leq i \leq n$, occur at $i \in I$, 
and the jumps of 
the sequence $\operatorname{dim}(X \cap F_i)$ occur at $i \in I_P$, so 
$X \in \Omega_{I_P}^{\circ}(F_{\sssize{\bullet}})$, which proves (i).  

For (ii), using Lemma 1(ii) similarly, we may assume $U \in  
\Omega_I^\circ(F_{\sssize{\bullet}})$ and $Y \in 
\Omega_P^\circ(F_{\sssize{\bullet}} W)$.  The jumps in the sequence 
$\operatorname{dim}(Y \cap F_kW)$, $1 \leq k \leq n-r$, occur at $k$ in 
$P$.  So the jumps in the sequence 
$$
\operatorname{dim}(Z \cap (F_i + U)), \quad 1 \leq i \leq n,
$$
occur at $i$ in $I_P^c$.   It follows that the jumps in the sequence 
$\operatorname{dim}(Z \cap F_i)$, $1 \leq i \leq n$, must occur at $i$ in 
$I \cup I_P^c$.  For if $i \notin I \cup I_P^c$, then $U \cap F_i = U \cap 
F_{i-1}$ and 
$Z \cap (F_i + U) = Z \cap (F_{i-1} + U)$, and these imply that $Z \cap F_i 
= Z \cap F_{i-1}$.  Therefore $Z$ is in 
$\Omega_{I_P^{+}}^{\circ}(F_{\sssize{\bullet}})$, which proves (ii).  \enddemo

\proclaim{Proposition 1}  Suppose $I(1), \ldots , I(m)$ are subsets of $[n]$ of 
cardinality $r$, with $\prod_{s=1}^m \omega_{I(s)} \neq 0$ in 
$H^*(\operatorname{Gr}(r,n))$.
\roster\widestnumber\item{(ii)}
\item"{(i)}"  If $P(1), \ldots , P(m)$, subsets 
of $[r]$ of cardinality $x$, have 
$\prod_{s=1}^m \omega_{P(s)} \neq 0$ in $H^*(\operatorname{Gr}(x,r))$, 
then $\prod_{s=1}^m \omega_{I(s)_{P(s)}} \neq 0$ in 
$H^*(\operatorname{Gr}(x,n))$.
\item"{(ii)}" If $P(1), \ldots , P(m)$, subsets of $[n-r]$ of cardinality $y$, 
have $\prod_{s=1}^m \omega_{P(s)} \neq 0$ in 
$H^*(\operatorname{Gr}(y,n-r))$, then  
$\prod_{s=1}^m \omega_{I(s)_{P(s)}^{+}} \neq 0$ in 
$H^*(\operatorname{Gr}(r+y,n))$.
\endroster
\endproclaim

\demo{Proof}  We use the fact (cf. [7], \S 4) 
 that for $m$ general flags $F_{\sssize{\bullet}}(1), 
\ldots, F_{\sssize{\bullet}}(m)$, the intersection $\cap_{s=1}^m 
\Omega_{I(s)} F_{\sssize{\bullet}}(s)$ is not empty if and only if the product 
$\prod_{s=1}^m \omega_{I(s)}$ of their classes is not zero.  
Moreover, for arbitrary flags 
$F_{\sssize{\bullet}}(s)$, if the product of the classes is not zero, then the 
intersection $\cap_{s=1}^m \Omega_{I(s)} F_{\sssize{\bullet}}(s)$ cannot be  
empty.  

Take general flags $F_{\sssize{\bullet}}(1), \ldots, 
 F_{\sssize{\bullet}}(m)$.   
There is an $r$-dimensional subspace $U$ that is in $\cap_{s=1}^m 
\Omega_{I(s)} F_{\sssize{\bullet}}(s)$.  In case (i), there is an $X$ in 
$\cap_{s=1}^m \Omega_{P(s)} F_{\sssize{\bullet}}(s) U$.  By Lemma 2(i), 
$X$ is in $\cap_{s=1}^m \Omega_{I(s)_{P(s)}} F_{\sssize{\bullet}}(s)$, 
which proves (i).  In case (ii), there is a $Y = Z/U$ in  $\cap_{s=1}^m 
\Omega_{P(s)} F_{\sssize{\bullet}}(s) W$, with $W = V/U$, and, by Lemma 
2(ii), $Z$ is in  $\cap_{s=1}^m \Omega_{I(s)_{P(s)}^{+}} 
F_{\sssize{\bullet}}(s)$; this proves (ii).   \enddemo

\remark{Remark}  In case (i), 
even if both products are the classes of a point, that is, 
if $\prod_{s=1}^m \omega_{I(s)} = \omega_{[r]}$ and $\prod_{s=1}^m 
\omega_{P(s)} = \omega_{[x]}$, it does not follow that $\prod_{s=1}^m 
\omega_{I(s)_{P(s)}} = \omega_{[x]}$, or even that the codimension of this 
class must be $x(n-x)$.  For example, with $n = 4$, $r = 2$, $x = 1$, and 
$\mathscr{I} = (\{2,4\},\{2,4\},\{2,3\})$, $\mathscr{P} = (\{2\},\{2\},\{1\})$, we have 
$\mathscr{I}_{\mathscr{P}} = (\{4\},\{4\},\{2\})$, whose product has 
codimension $2$, not $3$.   Similarly in case (ii), one can take $\mathscr{I} 
= (\{2,4\},\{2,4\},\{1,4\})$, $\mathscr{P} = (\{2\},\{2\},\{1\})$, with 
$\mathscr{I}_{\mathscr{P}}^{+} = (\{2,3,4\},\{2,3,4\},\{1,2,4\})$, 
again with the product of codimension $2$, not $3$.
\endremark
\smallpagebreak

We will use standard notation and facts about Littlewood-Richardson 
coefficients $c_{\alpha\,\beta}^{\,\,\gamma}$, for which the discussion in 
[7], \S 3 should suffice.  Note that if the lengths of the three partitions are at 
most $n$, and their widths (their first entries) are at most $N-n$, then 
$c_{\alpha\,\beta}^{\,\,\gamma}$ is the coefficient of the class 
$\sigma_{\gamma}$ in the product $\sigma_{\alpha} \cdot 
\sigma_{\beta}$ in $H^*(\operatorname{Gr}(n,N))$.  We will also need the 
following lemma. 

\proclaim{Lemma 3}  Let $\alpha$, $\beta$, $\gamma$ be partitions, with 
$c_{\alpha\,\beta}^{\,\,\gamma} > 0$.  Let $\widetilde{\gamma}$ be a 
partition with $\widetilde{\gamma} \subset \gamma$.  Then there are 
partitions $\widetilde{\alpha} \subset \alpha$ and $\widetilde{\beta} 
\subset \beta$ with $c_{\widetilde{\alpha} \, \widetilde{\beta}}^{\,\, 
\widetilde{\gamma}} > 0$.  
\endproclaim

\demo{Proof}\!\!\footnote{We did not find this fact in the literature, 
although, 
when asked, Buch, S. Fomin, J. Stembridge, Tao, and A. Zelevinski 
quickly produced different proofs of a more combinatorial flavor.  Buch 
shows in fact that if $\alpha \subset \widetilde{\gamma}$, then one can find 
$\widetilde{\beta} \subset \beta$ with  $c_{\alpha \,\widetilde{\beta}}^{\,\, 
\widetilde{\gamma}} > 0$.}   We deduce this from a theorem of Green and 
Klein [9], which says that for some (or any) discrete valuation ring, there is 
a finitely generated torsion module $C$ of type $\gamma$, with a submodule 
$B$ of type $\beta$, whose quotient module $A = C/B$ has type $\alpha$, if 
and only if $c_{\alpha\,\beta}^{\,\,\gamma}$ is positive.  
Choose a submodule 
$\widetilde{C}$ of $C$ of type $\widetilde{\gamma}$.  Then $\widetilde{B} 
= \ B \cap \widetilde{C}$ has type $\widetilde{\beta} \subset \beta$, and, 
since $\widetilde{A} = \widetilde{C}/\widetilde{B} \hookrightarrow A$, 
$\widetilde{A}$ 
has type $\widetilde{\alpha} \subset \alpha$.  By the Green-Klein theorem 
again, $c_{\widetilde{\alpha} \,\widetilde{\beta}}^{\, \widetilde{\gamma}} > 
0$.  \enddemo

\head Section 2.  General theorems and proofs \endhead

In the statements of theorems in this section, phrases in brackets indicate 
alternative versions, meaning that the theorem is true with or without any or 
all of these bracketed phrases.

For $1 \leq r \leq n$ and $m \geq 1$, as in [7] we let $S_r^n(m)$ be 
the set of $m$-tuples $\mathscr{I} = (I(1), \ldots , I(m))$ of subsets of 
cardinality $r$ in $[n] = \{ 1, \ldots , n \}$ such that the product of the 
corresponding classes $\omega_{I(s)}$ in $H^*(\operatorname{Gr}(r,n))$ 
does not vanish:
$$
S_r^n(m)\, =\, \bigl\{ \mathscr{I} 
= (I(1), \ldots , I(m)) \, \mid \, \prod_{s=1}^m \omega_{I(s)} \neq 0 \bigr\}. 
$$
Let $R_r^n(m)$ be the set of  $m$-tuples $\mathscr{I} $ whose product is 
the class of a point in the Grassmannian, with coefficient $1$: 
$$
R_r^n(m) \,=\, 
\bigl\{ \mathscr{I} = (I(1), \ldots , I(m)) \, \mid \, \prod_{s=1}^m \omega_{I(s)} = 
\omega_{[r]} \in H^{2r(n-r)}(\operatorname{Gr}(r,n))\bigr\}.
$$
We augment these lists for $r = n$, so that $R_n^n(m) = S_n^n(m)$ consists 
of the set $[n]$ repeated $m$ times.  We set
$$
S^n(m) = \bigcup_{1 \leq r \leq n} S_r^n(m), \quad {\text {and}} \quad
R^n(m) = \bigcup_{1 \leq r \leq n} R_r^n(m).
$$

The symmetric version of Theorem 1, for any number of factors, is:

\proclaim{Theorem 3}  Let $\alpha(1), \ldots , \alpha(m)$ be sequences of 
$n$-tuples of real numbers, with $\alpha(s) = (\alpha_1(s) \geq \alpha_2(s) 
\geq \ldots \geq \alpha_n(s))$.  There are complex Hermitian [real 
symmetric] [quaternionic Hermitian] $n$ by $n$ matrices $A(1), \ldots , 
A(m)$ with eigenvalues $\alpha(1), \ldots , \alpha(m)$ and $A(1) + \ldots + 
A(m) \leq 0$ if and only if 
$$
\sum_{s=1}^m \sum_{i \in I(s)} \alpha_i(s) \leq 0  \tag 
{$\dagger_{\mathscr{I}}$}
$$
for all $\mathscr{I} = (I(1), \ldots , I(m))$ in $S^n(m)$  
[$R^n(m)$].  \endproclaim

\demo{Proof}  The necessity of the conditions follows immediately from 
Klyachko's theorem. 
\footnote{It also follows almost as easily from the proof of the easy half of 
this theorem: Choose flags $F_{\sssize{\bullet}}(s)$ with $F_k(s)$ spanned 
by eigenvectors of $A(s)$ corresponding to its first $k$ eigenvalues.  Since 
$\prod_{s=1}^m \omega_{I(s)} \neq 0$, there is an $L$ in $\cap_{s=1}^m 
\Omega_{I(s)}(F_{\sssize{\bullet}}(s))$.  If $R_A(L)$ denotes the Rayleigh 
trace (the trace of $L \to \mathbb{C}^n \to \mathbb{C}^n \to L$, where the 
first map is the inclusion, the second is given by $A$, and the third is 
orthogonal projection onto $L$), then by the easy Hersch-Zwahlen Lemma 
(cf. [7], Prop. 1),  with $A = \sum_{s=1}^m A(s)$,
$$
\sum_{s=1}^m \sum_{i \in I(s)} \alpha_i(s) \leq \sum_{s=1}^m R_{A(s)}(L) 
= R_A(L) \leq 0,
$$
the last inequality since $A$ is negative semidefinite.}
Indeed, if $\sum_{s=1}^m A(s) \leq 0$, there is a positive semidefinite 
$A(m+1)$ with $\sum_{s=1}^{m+1} A(s) = 0$.  For any $\mathscr{I} = (I(1), 
\ldots , I(m))$ in $S_r^n(m)$, the $(m+1)$-tuple $\mathscr{I}' = (I(1), \ldots , I(m),\{ n-r+1, 
\ldots , n \})$ is in $S_r^{n}(m+1)$ (since $\omega_{\{ n-r+1, \ldots , n 
\}} = 1$), and the inequality ($\dagger_{\mathscr{I}'}$) implies 
($\dagger_{\mathscr{I}}$) since all the eigenvalues of $A(m+1)$ are 
nonnegative.  Note that for $r = n$, the inequality says that $\sum_{s=1}^m 
\operatorname{tr}(A(s)) = \operatorname{tr} (\sum_{s=1}^m A(s)) \leq 0$. 

We prove the converse by induction on $n$, the case $n = 1$ being trivial.  
Assume $\alpha(1), \ldots , \alpha(m)$ satisfy the inequalites 
($\dagger_{\mathscr{I}}$) for all $\mathscr{I} \in S^n(m)$.  We consider 
first the case where there is some $\mathscr{I} = (I(1), \ldots , I(m))$ in 
some $S_r^n(m)$ for some $1 \leq r \leq n$ for which the inequality 
($\dagger_{\mathscr{I}}$) is satisfied with equality: $\sum_{s=1}^m 
\sum_{i \in I(s)} \alpha_i(s) = 0$.  Define $r$-tuples $\alpha'(s)$, $1 \leq s 
\leq m$, by taking the $r$ values $\alpha_i(s)$ for $i \in I(s)$.  In 
symbols, if $I(s) = \{ i_1(s) <
\ldots < i_r(s) \}$, then $\alpha'(s)$ is
$$
\alpha_{i_1(s)} \geq \ldots \geq \alpha_{i_r(s)}. 
$$
Define similarly $\alpha''(s)$, $1 \leq s \leq m$, by taking the $n-r$ values 
$\alpha_i(s)$ for $i \notin I(s)$.  

\proclaim{Claim}  The $m$-tuple $\alpha'(1), \ldots , \alpha'(m)$ satisfies 
the inequalites $(\dagger_{\mathscr{P}})$ for all $\mathscr{P} \in 
S^r(m)$, and the  $m$-tuple $\alpha''(1), \ldots , \alpha''(m)$ satisfies the 
inequalites $(\dagger_{\mathscr{P}})$ for all $\mathscr{P} \in S^{n-r}(m)$.
\endproclaim

To prove the claim, for any $\mathscr{P} \in S^r(m)$, Proposition 1(i) 
implies that the $m$-tuple $\mathscr{I}_{\mathscr{P}}$ is in $S^r(n)$.  
The inequality ($\dagger_{\mathscr{I}_{\mathscr{P}}}$) for $\alpha(1), 
\ldots , \alpha(m)$ is exactly the inequality ($\dagger_{\mathscr{P}}$) for  
$\alpha'(1), \ldots , \alpha'(m)$.  Similarly by Proposition 1(ii), 
$\mathscr{I}_{\mathscr{P}}^{+}$ is in $S^n(m)$, and the inequality 
($\dagger_{\mathscr{I}_{\mathscr{P}}^{+}}$) for $\alpha(1), \ldots , 
\alpha(m)$ is equivalent to the inequality ($\dagger_{\mathscr{P}}$) for  
$\alpha''(1), \ldots , \alpha''(m)$, as one sees by adding the number 
$\sum_{s=1}^m\sum_{i \in I(s)} \alpha_i(s) = 0$ to the left side.
\smallpagebreak

By the claim, the $m$-tuple $\alpha'(1), \ldots , \alpha'(m)$ satisfies all the 
inequalities ($\dagger_{\mathscr{P}}$) for $\mathscr{P} \in S^r(m)$ of 
Klyachko's theorem. By the real version of this theorem ([6] Thm. 4, cf 
[7] \S 10.7), there are real symmetric $r$ by $r$ matrices $A'(1), \dots ,
A'(m)$ with eigenvalues $\alpha'(1), \ldots , \alpha'(m)$ and 
$\sum_{s=1}^m A'(s) = 0$.  Similarly, the $\alpha''(1), \ldots , \alpha''(m)$ 
satisfy all the inequalities ($\dagger_{\mathscr{P}}$) for $\mathscr{P} \in 
S^{n-r}(m)$.  By induction, since $n-r < n$, there are $n-r$ by $n-r$ real 
symmetric matrices $A''(1), \dots, A''(m)$ with eigenvalues $\alpha''(1), 
\ldots , \alpha''(m)$ and $\sum_{s=1}^m A''(s) \leq 0$.  The direct sum 
matrices $A(s) = \left(\smallmatrix A'(s) & 0 \\ 0 &A''(s) \endsmallmatrix 
\right)$ are then real symmetric $n$ by $n$ matrices whose eigenvalues 
are the given 
$\alpha(1), \ldots , \alpha(m)$, with $A(1) + \ldots + A(m) \leq 0$.

Now suppose that all the inequalities ($\dagger_{\mathscr{I}}$) are strict, 
for all $\mathscr{I} \in S^n(m)$.  Since this is a finite list, we may find a 
positive $\epsilon$ such that, after replacing each $\alpha_i(s)$ by 
$\widetilde{\alpha}_i(s) = \alpha_i(s) + \epsilon$, all the inequalities 
($\dagger_{\mathscr{I}}$) are satisfied for  $\widetilde{\alpha}(1), \ldots , 
\widetilde{\alpha}(m)$, but at least one holds with equality.  By the case just 
proved, there are real symmetric $\widetilde{A}(1), \ldots ,
\widetilde{A}(m)$ with eigenvalues $\widetilde{\alpha}(1), \ldots , 
\widetilde{\alpha}(m)$ and $\sum_{s=1}^m\widetilde{A}(s) \leq 0$.  The 
matrices $A(s) = \widetilde{A}(s) - D(\epsilon)$, where $D(\epsilon)$ is the 
diagonal matrix with diagonal entries $\epsilon$, satisfy the required 
conditions.  

To complete the proof of the theorem, we must show that the inequalites 
($\dagger_{\mathscr{I}}$) for $\mathscr{I} \in S^n(m)$ follow from those for 
$\mathscr{I} \in R^n(m)$.  This is proved by the method of Belkale and 
Woodward.  In fact, the proof of Proposition 10 in [7], \S 7 shows that, 
for any weakly decreasing sequences  $\alpha(1), \ldots , \alpha(m)$ of $n$ 
real numbers, 
$$
\max_{\mathscr{I} \in S^n(m)} \sum_{s=1}^m\sum_{i \in I(s)} \alpha_i(s) 
= \max_{\mathscr{I} \in R^n(m)} \sum_{s=1}^m\sum_{i \in I(s)} 
\alpha_i(s).
$$
\enddemo
\proclaim{Corollary}  Given 
weakly decreasing $n$-tuples $\alpha(1), \ldots , \alpha(m)$,  
$\gamma$, there are complex Hermitian [real symmetric] [quaternionic 
Hermitian] $n$ by $n$ matrices $A(1),\! \ldots \! , A(m)$, $C$  with 
eigenvalues $\alpha(1), \ldots , \alpha(m)$, $\gamma$ with $C \leq 
A(1) + \ldots + A(m)$ if and only if 
$$
\sum_{k \in K} \gamma_k \, \leq \, \sum_{s=1}^m \sum_{i \in I(s)} \alpha_i(s) 
\tag{$\dagger_{\mathscr{I} K}$}
$$
for all $I(1), \ldots , I(m), K$ of cardinality $r$ in $[n]$ such that 
$\sigma_{\lambda(K)}$ occurs in the product $\prod_{s=1}^m 
\sigma_{\lambda(I(s))}$ [with coefficient $1$] in 
$H^*(\operatorname{Gr}(r,n))$, for all $r \leq n$.
\endproclaim
\demo{Proof}   The theorem for $m+1$ factors is applied to the situation 
$$
- A(1) - A(2) - \ldots - A(m) + C \leq 0
$$
which changes signs and reverses the order of the eigenvalues of each 
$A(s)$.  One uses duality in the Schubert calculus ($\sigma_{\lambda(I)}$ is 
dual to $\omega_I$) to make the translation. For details on this translation, see [7], \S 10.1. Note that the inequality 
$\sum_{i=1}^n \gamma_i \leq \sum_{s=1}^m \sum_{i=1}^n \alpha_i(s)$ is that 
for $I(1) = \ldots = I(m) = K = [n]$, with $r = n$. 
\enddemo

Theorem 1 from the introduction is a special case of this corollary.

If $A$ is a compact, self-adjoint, and positive semidefinite, linear operator 
on a separable Hilbert space, there is an orthonormal basis $e_1, e_2, 
\ldots $ so that $A\,e_i = \alpha_i e_i$ for all $i$, with $\alpha_1 \geq 
\alpha_2 \geq \ldots$, and $\lim \alpha_i = 0$; these $\alpha_i$ are the  
eigenvalues of $A$.  The operator is said to be {\em of trace class} if 
$\sum_{i=1}^{\infty} \alpha_i < \infty$.  The proof given by Friedland in 
[4], \S 5 applies without change to give the following stronger result: 

\proclaim{Theorem 4}  Let  $\alpha(1), \ldots , \alpha(m)$ and $\gamma$ 
be weakly decreasing infinite sequences of nonnegative numbers that 
converge to zero. 
These occur as eigenvalues of compact self-adjoint operators 
$A(1), \ldots , A(m)$ and $C$ on a real or complex separable Hilbert 
space, with $C \leq A(1) + \ldots + A(m)$, if and only if they satisfy the 
inequalities $(\dagger_{\mathscr{I} K})$ 
for all $I(1), \ldots , I(m), K$ of cardinality $r$ in $[n]$ such that 
$\sigma_{\lambda(K)}$ occurs in the product $\prod_{s=1}^m 
\sigma_{\lambda(I(s))}$ [with coefficient $1$] in 
$H^*(\operatorname{Gr}(r,n))$, for all $1 \leq r \leq n < \infty$. 
Moreover, if each sum $\sum_{i=1}^{\infty} \alpha_i(s) $ is finite 
and $\sum_{i=1}^{\infty} \gamma_i = \sum_{s=1}^m \sum_{i=1}^{\infty}
\alpha_i(s)$, then these inequalities are necessary and sufficient for 
the existence of operators $A(1), \ldots, A(m)$ and $C$ 
with these eigenvalues, and $C = A(1) + \ldots + A(m)$.
 \endproclaim

Note that an inequality 
$(\dagger_{\mathscr{I} K})$ depends only on $\mathscr{I}$ and $K$, and 
not on $n$; to get a list without repeats, one may take $n$ to be the 
largest integer in $K$.  
In contrast to the case for matrices, as Friedland points out, the 
resulting list is far from minimal.  Indeed, any finite number of the 
inequalities $(\dagger_{\mathscr{I} K})$ can be omitted.  For $m=2$, 
the fact that $(\dagger_{I J K})$ is redundant follows from 
the fact that for $(I,J,K) \in S_r^n$, and any $a > n$ and $b > n$, the 
triple $(I\cup\{a\},J\cup\{b\},K\cup\{c\})$ is in $S_{r+1}^c$, where 
$c = a+b-r-1$.  Indeed, the equality
$$
c_{\lambda(I\cup\{a\}) \, \lambda(J\cup\{b\})}^{\,\,\,\,\,\,\,\,
\lambda(K\cup\{c\})} =
c_{\lambda(I) \, \lambda(J)}^{\,\,\,\,\lambda(K)}
$$
follows from a simple general fact about Littlewood-Richardson coefficients: 
if $\lambda$, $\mu$, and $\nu$ are partitions, with $\lambda_1 + \mu_1 = 
\nu_1$, then $c_{\lambda \, \mu}^{\,\,\nu}$ is equal to the 
Littlewood-Richardson coefficient for the three partitions obtained 
from $\lambda$, $\mu$, and $\nu$ by removing each of their first 
entries (cf. [5], Chapter 5).  Letting $a$, $b$, and therefore $c$ 
tend to infinity, one sees that $(\dagger_{I J K})$ follows from the 
inequalities $(\dagger_{I\cup\{a\}\, J\cup\{b\}\, K\cup\{c\}})$.

\proclaim{Proposition 2}  Let $\alpha(1), \ldots , \alpha(m)$ be weakly 
decreasing sequences of $n$ integers.  The following are equivalent:
\roster
\item"{(1)}" The inequality $(\dagger_{\mathscr{I}})$ is valid for all 
$\mathscr{I}$ in $S^n(m)$ [$R^n(m)$].
\item"{(2)}" There are weakly decreasing sequences $\widetilde{\alpha}(1), 
\ldots , \widetilde{\alpha}(m)$ of integers such that $\widetilde{\alpha}_i(s) 
\geq \alpha_i(s)$ for $1 \leq s \leq m$, $1 \leq i \leq n$, with $\sum_{s=1}^m 
\sum_{i=1}^n \widetilde{\alpha}_i(s) = 0$, such that 
$\widetilde{\alpha}(1), \ldots , \widetilde{\alpha}(m)$ satisfies 
$(\dagger_{\mathscr{I}})$ for all $\mathscr{I}$ in $S^n(m)$ [$R^n(m)$].
\item"{(3)}" The same assertion as in $(2)$, except that 
$\widetilde{\alpha}(s) = \alpha(s)$ 
for all $s \neq s_0$, for any choice of $s_0 
\in \{ 1, \ldots , m \}$.
\endroster
\endproclaim
\demo{Proof}  
The assertions (3) $\Rightarrow$ (2) $\Rightarrow$ (1) are clear.  
To prove that (1) $\Rightarrow$ (3), it suffices to show that if $\sum_{s=1}^m 
\sum_{i \in I(s)} \alpha_i(s) < 0$, we can increase one of the integers 
$\alpha_i(s_0)$ by $1$, leaving the others unchanged, so that all the 
inequalities ($\dagger_{\mathscr{I}}$) remain valid; repeating this until the 
full sum becomes $0$ produces the required $\widetilde{\alpha}(s_0)$.  If 
all the ($\dagger_{\mathscr{I}}$) are strict inequalities, we can simply 
replace $\alpha_1(s_0)$ by  $\alpha_1(s_0) + 1$.  Otherwise take $r$ 
maximal such that ($\dagger_{\mathscr{I}}$) is an equality for some 
$\mathscr{I} \in S_r^n(m)$.  As in the proof of Theorem 3, this partitions 
the $m$-tuple $\boldsymbol{\alpha} = ({\alpha}(1), \ldots , 
{\alpha}(m))$ into two $m$-tuples $\boldsymbol{\alpha}'$ and 
$\boldsymbol{\alpha}''$, satisfying the 
conditions ($\dagger_{\mathscr{P}}$) 
for $\mathscr{P}$ in $S^r(m)$ and $S^{n-r}(m)$ respectively, and with 
$\sum_{s=1}^m \sum_{i=1}^r \alpha_i^{\prime}(s) = 0$.  By the maximality 
of $r$, all the inequalities ($\dagger_{\mathscr{P}}$) for 
$\boldsymbol{\alpha}''$ 
must be strict; otherwise ($\dagger_{\mathscr{I}_{\mathscr{P}}^{+}}$) 
would be an equality for 
$\boldsymbol{\alpha}$, contradicting the maximality of $r$.

Let ${\widetilde{\boldsymbol{\alpha}}''}$ 
be obtained from $\boldsymbol{\alpha}''$ by 
increasing $\alpha_1^{\prime \prime}(s_0)$ by $1$, leaving all other 
values unchanged.  By Theorem 3, $\boldsymbol{\alpha}'$ and 
${\widetilde{\boldsymbol{\alpha}}''}$ are 
eigenvalues of Hermitian matrices $A'(1), 
\ldots , A'(m)$ and $\widetilde{A}''(1), \ldots , \widetilde{A}''(m)$, with 
$\sum_{s=1}^m A'(s) = 0$ and $\sum_{s=1}^m \widetilde{A}''(s) \leq 0$.  
The direct sums $A(s) = A'(s) \oplus \widetilde{A}''(s)$ also have 
$\sum_{s=1}^m A(s) \leq 0$, so by Theorem 3 their eigenvalues satisfy the 
inequalities ($\dagger_{\mathscr{I}}$) for all $\mathscr{I} \in S^n(m)$.  
The eigenvalues of $A(s)$ are $\alpha(s)$ for $s \neq s_0$, while the 
eigenvalues of $A(s_0)$ are obtained from $\alpha(s_0)$ by increasing one 
$\alpha_i(s_0)$ by $1$  (this $i$ is the smallest integer 
in $[n]$ not in $I(s_0)$).  
\enddemo

\proclaim{Theorem 5}  Let $\alpha(1), \ldots , \alpha(m)$ and $\gamma$ 
be partitions of lengths at most $n$.  The following are equivalent:
\roster
\item"{(1)}"  The inequalities $(\dagger_{\mathscr{I} K})$
are valid for all subsets $I(1), \ldots , I(m), K$ of $[n]$ of cardinality $r$, $1 
\leq r \leq n$, such that $\sigma_{\lambda(K)}$ occurs in $\prod_{s=1}^m 
\sigma_{\lambda(I(s)}$ [with coefficient $1$] in 
$H^*(\operatorname{Gr}(r,n))$.
\item"{(2)}" There are complex Hermitian [real symmetric] [quaternionic 
Hermitian] matrices $A(1), \ldots , A(m)$ and $C$ with eigenvalues 
$\alpha(1), \ldots , \alpha(m)$ and $\gamma$ such that $C \leq A(1) + 
\ldots + A(m)$.
\item"{(3)}" There is a partition $\widetilde{\gamma} \supset \gamma$ 
such that $\sum_{i=1}^n \widetilde{\gamma}_i = \sum_{s=1}^m 
\sum_{i=1}^n \alpha_i(s)$ and the inequalities $(\dagger_{\mathscr{I} 
K})$  hold for the same $\mathscr{I}$, $K$ as in $(1)$ for all $r < n$.
\item"{(4})" There are partitions $\widetilde{\alpha}(s) \subset \alpha(s)$ 
so that $\sum_{i=1}^n \gamma_i = \sum_{s=1}^m \sum_{i=1}^n 
\widetilde{\alpha}_i(s)$ and the inequalities $(\dagger_{\mathscr{I} K})$  
hold for the same $\mathscr{I}$, $K$ as in $(1)$ for all $r < n$. 
\endroster
\endproclaim
\remark{Remark} The partitions $\alpha(1), \ldots , \alpha(m)$ and 
$\widetilde{\gamma}$ in (3), as well as the partitions $\widetilde{\alpha}(1), 
\ldots , \widetilde{\alpha}(m)$ and $\gamma$ of (4), satisfy the Klyachko 
conditions, so all the equivalent conditions of [7], \S 10, Thm\. 17 can be 
substituted, and the list continued.  For example, with $V(\lambda)$ 
denoting the irreducible (polynomial) representation of 
$\operatorname{GL}(n,\mathbb{C})$ of highest weight $\lambda$, 
condition (3) and (4) are equivalent to:
\proclaim\nofrills{}\roster
\item"{(5)}"  There is a partition $\widetilde{\gamma} \supset \gamma$ 
such that the representation $V(\widetilde{\gamma})$ occurs in 
$V(\alpha(1))\otimes \ldots \otimes V(\alpha(m))$.
\item"{(6)}"  There are partitions $\widetilde{\alpha}(s) \subset \alpha(s)$ 
such that the representation 
$V(\gamma)$ occurs in  $V(\widetilde{\alpha}(1))\otimes \ldots 
\otimes V(\widetilde{\alpha}(m))$.  \endroster \endproclaim
There are also equivalent conditions involving invariant factors: 
\proclaim\nofrills{}\roster
\item"{(7)}"  For some [every] discrete valuation ring $R$, there exist $n$ by 
$n$ matrices $A(1), \ldots , A(m)$ with entries in $R$ and invariant 
factors $\alpha(1), \ldots , \alpha(m)$ such that the product $A(1) \cdot 
\ldots \cdot A(m)$ has invariant factors $\widetilde{\gamma}$ for some 
$\widetilde{\gamma} \supset \gamma$.
\item"{(8)}"  For some [every] discrete valuation ring $R$, there is an 
$R$-module 
$C$ of type $\gamma$, with a filtration $0 = C(0) \subset C(1) 
\subset \ldots \subset C(m) = C$ of submodules such that $C(s)/C(s-1)$ is 
isomorphic to a submodule of a module of type $\alpha(s)$, for $1 \leq s 
\leq m$.  \endroster \endproclaim 
\endremark
   
\demo{Proof}  The equivalence of (1) and (2) is a case of the Corollary to 
Theorem 3.  Similarly, the equivalence of (1) and (3) of Proposition 2, 
applied with $m+1$ factors and with $s_0 = m+1$ corresponding to 
$\gamma$, shows that (1) and (3) are equivalent.  Since (4) clearly implies 
(1), it remains to show that (3) implies (4).  This amounts to the following 
assertion.  We are given a partition $\widetilde{\gamma} \supset \gamma$ 
such that $\sigma_{\widetilde{\gamma}}$ appears in the product 
$\prod_{s=1}^m \sigma_{\alpha(s)}$ in $H^*(\operatorname{Gr}(n,N))$, 
for any $N$ with $N-n \geq \widetilde{\gamma}_1$.  We must find 
partitions $\widetilde{\alpha}(s) \subset \alpha(s)$ so that 
$\sigma_{\gamma}$ appears in $\prod_{s=1}^m 
\sigma_{\widetilde{\alpha}(s)}$.  We argue by induction on $m$, the case 
$m = 2$ being Lemma 3 in Section 1.  For $m > 2$, there is some 
$\sigma_{\beta}$ occuring in $\prod_{s > 1} \sigma_{\alpha(s)}$ such that 
$\sigma_{\widetilde{\gamma}}$ appears in $\sigma_{\alpha(1)} \cdot 
\sigma_{\beta}$.  By the case for $m = 2$, we can find 
$\widetilde{\alpha}(1) \subset \alpha(1)$ and $\widetilde{\beta} \subset 
\beta$ so that $\sigma_{\gamma}$ occurs in 
$\sigma_{\widetilde{\alpha}(1)} \cdot \sigma_{\widetilde{\beta}}$.  By the 
case for $m-1$, there are $\widetilde{\alpha}(s) \subset \alpha(s)$, $2 \leq 
s \leq m$, so that $\sigma_{\widetilde{\beta}}$ occurs in $\prod_{s > 1} 
\sigma_{\widetilde{\alpha}(s)}$.  Therefore $\sigma_{\gamma}$ occurs in 
$\prod_{s=1}^m \sigma_{\widetilde{\alpha}(s)}$, as required.  \enddemo

Theorem 2 of the introduction follows from this theorem.  Indeed $\alpha$, 
$\beta$, $\gamma$ satisfy the conditions of Theorem 2 exactly when there 
are $\widetilde{\alpha} \subset \alpha$ and $\widetilde{\beta} \subset 
\beta$ such that the Littlewood-Richardson coefficient 
$c_{\widetilde{\alpha}\, \widetilde{\beta}}^{\,\,\gamma}$ is positive.  By 
the Green-Klein theorem [9] this is equivalent to the existence of a short 
exact sequence of modules $0 \to \widetilde{B} \to C \to \widetilde{A} \to 
0$, with $\widetilde{A}$ of type $\widetilde{\alpha}$,  $\widetilde{B}$ of 
type $\widetilde{\beta}$, and $C$ of type $\gamma$.  Taking $A$ of type 
$\alpha$ and $B$ of type $\beta$, there is an epimorphism $B \to 
\widetilde{B}$ and a monomorphism $\widetilde{A} \to A$, and this gives 
the required exact sequence $B \to C \to A$.
\smallpagebreak
        
We now turn to the question of the minimality of these lists of inequalities.  
Knutson, Tao, and Woodward [12] have announced that, for $n \geq 3$, 
the inequalities ($\dagger$) and ($\dagger_{IJK}$), for $(I,J,K) \in R_r^n$ 
and $r < n$, when restricted to the hyperplane defined by $\sum_{i=1}^n 
\gamma_i = \sum_{i=1}^n \alpha_i + \sum_{i=1}^n \beta_i$, are 
independent: if any one is left out, the cone defined by the others is strictly 
larger.  It is not hard to deduce from this the independence of our 
inequalities, where the hyperplane condition is replaced by an inequality to 
be on one side of it:

\proclaim{Theorem 6}  For any $n \geq 1$, the inequalities $(\dagger)$ 
and $(\dagger_{IJK})$ for $(I,J,K) \in R_r^n$ and $1 \leq r \leq n$, are 
independent.  \endproclaim

\demo{Proof}  For $n = 1$ we have only the inequality $\gamma_1 \leq 
\alpha_1 + \beta_1$.  For $n = 2$ there are seven inequalities $\alpha_1 
\geq \alpha_2$, $\beta_1 \geq \beta_2$, $\gamma_1 \geq \gamma_2$, 
$\gamma_1 \leq \alpha_1 + \beta_1$, $\gamma_2 \leq \alpha_2 + 
\beta_1$, $\gamma_2 \leq \alpha_1 + \beta_2$, and $\gamma_1 + 
\gamma_2 \leq \alpha_1 + \alpha_2 + \beta_1 + \beta_2$, which are 
easily verified to be independent.  For $n \geq 3$, we know that none of the 
chamber inequalities ($\dagger$) or the inequalities ($\dagger_{IJK}$) for 
$(I,J,K) \in R_r^n$ and $r < n$ can be omitted; indeed, if one is omitted, 
the others, even when restricted to the hyperplane $\sum_{i=1}^n 
\gamma_i = \sum_{i=1}^n \alpha_i + \sum_{i=1}^n \beta_i$, define a 
larger cone.  Finally, the last inequality ($\dagger_{IJK}$), for $I = J = K = 
[n]$, cannot be omitted, since for example the triple $\alpha = \beta = 
\gamma = (n-1, n-3, \ldots , -n+3, -n+1)$ is in the boundary hyperplane but 
all the other inequalities are strict (cf. [6], Lemma 2).  
\enddemo

Granting the corresponding assertion from [12], the inequalities 
($\dagger$) and ($\dagger_{\mathscr{I}}$) for $\mathscr{I} \in R^n(m)$ 
in Theorem 3 are independent.  Again the special case $n=2$, with its 
$2m+1$ inequalites, is handled directly.   The same holds for the 
inequalities ($\dagger$) and ($\dagger_{\mathscr{I} K}$)  (for those 
$\mathscr{I},K$ such that $\sigma_{\lambda(K)}$ appears in 
$\prod_{s=1}^m \sigma_{\lambda(I(s))}$) in Theorem 5. 

There are similar results when the majorization is taken from the other 
side.  Since no new ideas are required for the proof (Lemma 3 is not 
even needed), we leave the proof, as well as the statements of variations 
as in the remark after Theorem 5, and the corresponding assertions about 
the minimality of the list of inequalities, to the reader.

\proclaim{Theorem 7}  Given weakly decreasing sequences $\alpha(1), 
\ldots , \alpha(m)$ and $\gamma$ of real numbers of length $n$, the 
following are equivalent:
\roster
\item"{(1)}" $\sum_{s=1}^m \sum_{i \in I(s)} \alpha_i(s) \leq \sum_{k \in 
K} \gamma_k$ for all $I(1), \ldots , I(m), K$ of cardinality $r \leq n$ such 
that $\omega_K$ appears [with coefficient $1$] in $\prod_{s=1}^m 
\omega_{I(s)}$.
\item"{(2)}"  There are Hermitian [real symmetric] [quaternionic Hermitian] 
$n$ by $n$ matrices $A(1), \ldots , A(m), C$ with $A(1) + \ldots + A(m) 
\leq C$. \endroster
If all $\alpha_i(s)$ and $\gamma_i$ are integers, these are equivalent to:
\roster
\item"{(3)}"  There are integral $\widetilde{\alpha}(s)$ with 
$\widetilde{\alpha}_i(s) \geq \alpha_i(s)$ for $1 \leq i \leq n$ and for all 
$s$ [with $\widetilde{\alpha}(s) = \alpha(s)$ for all $s \neq s_0$, for any 
given $s_0$], such that $V(\gamma)$ occurs in 
$V(\widetilde{\alpha}(1))\otimes \ldots \otimes V(\widetilde{\alpha}(m))$.
\endroster
If all $\alpha(s)$ and $\gamma$ are partitions, these are equivalent to:
 \roster
\item"{(4)}"  There is a partition $\widetilde{\gamma} \subset \gamma$ 
such that $V(\widetilde{\gamma})$ occurs in $V(\alpha(1))\otimes \ldots 
\otimes V(\alpha(m))$.
\endroster
\endproclaim

\enddocument